\documentclass[twoside,10pt]{article}
\usepackage{amssymb}

\usepackage[pagebackref,colorlinks=true]{hyperref}  

\begin{document}

\setlength{\textwidth}{126mm} \setlength{\textheight}{180mm}
\setlength{\parindent}{0mm} \setlength{\parskip}{2pt plus 2pt}

\frenchspacing

\pagestyle{myheadings}

\markboth{Kostadin Gribachev, Mancho Manev}{Almost hypercomplex
pseudo-Hermitian manifolds and an example}

\date{2007/10/23}

\newtheorem{thm}{Theorem}[section]
\newtheorem{lem}[thm]{Lemma}
\newtheorem{prop}[thm]{Proposition}
\newtheorem{cor}[thm]{Corollary}
\newtheorem{probl}[thm]{Problem}

\newtheorem{defn}{Definition}[section]
\newtheorem{rem}{Remark}[section]
\newtheorem{exa}{Example}



\newcommand{\X}{\mathfrak{X}}
\newcommand{\B}{\mathcal{B}}
\newcommand{\s}{\mathfrak{S}}
\newcommand{\g}{\mathfrak{g}}
\newcommand{\W}{\mathcal{W}}
\newcommand{\Lgr}{\mathrm{L}}
\newcommand{\dd}{\mathrm{d}}

\newcommand{\pd}{\partial}
\newcommand{\ddx}{\frac{\pd}{\pd x^i}}
\newcommand{\ddy}{\frac{\pd}{\pd y^i}}
\newcommand{\ddu}{\frac{\pd}{\pd u^i}}
\newcommand{\ddv}{\frac{\pd}{\pd v^i}}

\newcommand{\diag}{\mathrm{diag}}
\newcommand{\End}{\mathrm{End}}
\newcommand{\im}{\mathrm{Im}}
\newcommand{\Id}{\mathrm{Id}}

\newcommand{\ie}{i.e.}
\newfont{\w}{msbm9 scaled\magstep1}
\def\R{\mbox{\w R}}
\newcommand{\norm}[1]{\left\Vert#1\right\Vert ^2}
\newcommand{\nN}{\norm{N}}
\newcommand{\ad}{{\rm ad}}

\newcommand{\nJ}[1]{\norm{\nabla J_{#1}}}
\newcommand{\thmref}[1]{Theorem~\eqref{#1}}
\newcommand{\propref}[1]{Proposition~\eqref{#1}}
\newcommand{\secref}[1]{\S\eqref{#1}}
\newcommand{\lemref}[1]{Lemma~\eqref{#1}}
\newcommand{\dfnref}[1]{Definition~\eqref{#1}}
\newcommand{\eqref}[1]{(\ref{#1})}

\frenchspacing

\hyphenation{Her-mi-ti-an ma-ni-fold}
\hyphenation{iso-tro-pic}


\title{Almost hypercomplex pseudo-Hermitian manifolds and a
4-dimensional Lie group \\ with such structure}

\author{Kostadin Gribachev, Mancho Manev\thanks{Corresponding
author}}

\maketitle

{\small
{\it Abstract.} Almost hypercomplex pseudo-Hermitian manifolds are
considered. Isotropic hyper-K\"ahler manifolds are introduced. A
4-parametric family of 4-dimensional manifolds of this type is
constructed on a Lie group. This family is characterized
geometrically. The condition a 4-manifold to be isotropic hyper-K\"ahler is given.\\

{\it Mathematics Subject Classification (2000):} 53C26, 53C15, 53C50, 53C55   \\
{\it Key words:} almost hypercomplex manifold, pseudo-Hermitian
metric, indefinite metric, Lie group }

\setcounter{tocdepth}{2} \tableofcontents


\section*{Introduction}

The general setting of this paper is inspired by the work of
D.~V.~Alekseevsky and S.~Marchiafava \cite{AlMa}. Our purpose is
to develop a parallel direction including indefinite metrics. More
precisely we combine the ordinary Hermitian metrics with the
so-called by us skew-Hermitian metrics with respect to the almost
hypercomplex structure.

In the first section we consider an appropriate decomposition of
the space of all bilinear forms on a vector space equipped with a
hypercomplex structure. Here we emphasize on a notion of the
skew-Hermitian metric. In fact, we construct three skew-Hermitian
metrics and one Hermitian, \ie{} a pseudo-Hermitian structure.

In the second we develop the notion of an almost hypercomplex
manifold with a pseudo-Hermitian structure and particularly the
so-called pseudo-hyper-K\"ahler\-ian and isotropic K\"ahler
structures.

Finally, in the third section we equip a 4-dimensional Lie group
with an almost hypercomplex pseudo-Hermitian structure and we
characterize it geometrically.


\section{Hypercomplex pseudo-Hermitian structures\\ on a vector space}

Let $V$ be a real $4n$-dimensional vector space. A (local) basis
on $V$ is denoted by $\left\{ \pd / \pd x^i, \pd / \pd y^i, \pd /
\pd u^i, \pd / \pd v^i \right\}$, $i=1,2$,$\dots,n$. Each vector
$\textbf{x}$ of $V$ is represented in the mentioned basis as
follows
   \begin{equation}\label{11}
\textbf{x}=x^i\ddx+y^i\ddy+u^i\ddu+v^i\ddv.
   \end{equation}

A standard hypercomplex structure on $V$ is defined as in
\cite{So}:
\begin{equation}\label{12}
\begin{array}{llll}
J_1\ddx=\ddy, &J_1\ddy=-\ddx, &J_1\ddu=-\ddv, &J_1\ddv=\ddu;
\\
J_2\ddx=\ddu, &J_2\ddy=\ddv, &J_2\ddu=-\ddx, &J_2\ddv=-\ddy;
\\
J_3\ddx=-\ddv, &J_3\ddy=\ddu, &J_3\ddu=-\ddy, &J_3\ddv=\ddx.
\end{array}
\end{equation}

The following properties about $J_i$ are direct consequences of
\eqref{12}
\begin{equation}\label{13}
\begin{array} {l}
J_1^2=J_2^2=J_3^2=-\Id,\\
J_1J_2=-J_2J_1=J_3,\quad J_2J_3=-J_3J_2=J_1,\quad
J_3J_1=-J_1J_3=J_2.
\end{array}
\end{equation}

If $x \in V$, \ie{}
$x(x^1,\dots,x^n;y^1,\dots,y^n;u^1,\dots,u^n;v^1,\dots,v^n)$ then
according to \eqref{12} and \eqref{13} we have
\[
\begin{array} {l}
J_1x(-y^1,\dots,-y^n;x^1,\dots,x^n;v^1,\dots,v^n;-u^1,\dots,-u^n),\\
J_2x(-u^1,\dots,-u^n;-v^1,\dots,-v^n;x^1,\dots,x^n;y^1,\dots,y^n),\\
J_3x(v^1,\dots,v^n;-u^1,\dots,-u^n;y^1,\dots,y^n;-x^1,\dots,-x^n).\\
\end{array}
\]

\begin{defn}[\cite{AlMa}]\label{d1}
1) A triple $H=(J_1,J_2,J_3)$ of anticommuting complex structures
on $V$ with $J_3=J_1J_2$ is called {\em a hypercomplex
structure\/} on $V$;

2) The 3-dimensional subspace $Q\equiv\langle H\rangle=\mathbb{R}
J_1+\mathbb{R} J_2+\mathbb{R} J_3$ of the space of endomorphisms
$\End V$ is called {\em a quaternionic structure\/} on $V$. It is
said that $H=(J_\alpha)$ is an {\em admissible basis\/} of $Q$.
\end{defn}
Note that two admissible bases $H$ and $H'$ of $Q=\langle
H\rangle=\langle H'\rangle$ are related by an orthogonal matrix in
SO(3).

The matrices of $J_1$ and $J_2$ are given in \cite{So} by
$(n\times n)$-sets of $(4\times 4)$-matrices
$\mathbf{J_\alpha}=\diag\left(I_\alpha, I_\alpha, \dots ,
I_\alpha\right)$, where $I_\alpha$ $(\alpha=1,2,3)$ are
respectively
\[
 I_1=\left(%
\begin{array}{cccc}
         0 & 1 & 0 & 0 \\
        -1 & 0 & 0 & 0 \\
         0 & 0 & 0 & -1 \\
         0 & 0 & 1 & 0
 \end{array}
 \right),\;
 I_2=\left(%
  \begin{array}{cccc}
        0 & 0 & 1 & 0 \\
         0 & 0 & 0 & 1 \\
        -1 & 0 & 0 & 0 \\
         0 & -1& 0 & 0
 \end{array}
 \right),\;\]
and consequently
\[
I_3=\left(%
 \begin{array}{cccc}
        0 & 0 &  0 & 1 \\
         0 & 0 & -1 & 0 \\
         0 & 1 &  0 & 0 \\
        -1 & 0 &  0 & 0
 \end{array}
 \right).
\]

The matrices $\mathbf{J_\alpha}$ of the complex structures
$J_\alpha$ $(\alpha=1,2,3)$ with respect to an admissible frame
for $H=(J_\alpha)$ are called {\em standard matrices\/}.

A bilinear form $f$ on $V$ is defined as ordinary, $f:\; V \times
V 
\rightarrow \mathbb{R}$. We denote by $\B(V)$ the set of all
bilinear forms on $V$. Each $f$ is a tensor of type $(0,2)$, and
$\B(V)$ is a vector space of dimension $16n^2$.

Let $J$ be a given complex structure on $V$. A bilinear form $f$
on $V$ is called {\em Hermitian\/} (respectively, {\em
skew-Hermitian\/}) with respect to $J$ if the identity
$f(Jx,Jy)=f(x,y)$ (respectively, $f(Jx,Jy)=-f(x,y)$ holds true.

\begin{defn}[\cite{AlMa}]\label{d2}
A bilinear form $f$ on $V$ is called {\em an Hermitian bilinear
form\/} with respect to  $H=(J_\alpha)$ if it is Hermitian with
respect to any complex structure $J_\alpha$, $\alpha=1,2,3$, \ie{}
\[
f(J_\alpha x,J_\alpha y)=f(x,y),\qquad \forall x, y \in V.
\]
\end{defn}

We will denote by $\B_\mathrm{H}(V)$ the set of all Hermitian
bilinear forms on $V$.

In \cite{GrMaDi} is introduced the notion of pseudo-Hermitian
bilinear forms, namely:
\begin{defn}[\cite{GrMaDi}]\label{d3}
A bilinear form $f$ on $V$ is called {\em a pseudo-Hermitian
bilinear form with respect to\/} $H=(J_1,J_2,J_3)$, if it is
Hermitian with respect to $J_\alpha$ and skew-Hermitian with
respect to $J_\beta$ and $J_\gamma$, \ie{}
\begin{equation}\label{16}
    f(J_\alpha x,J_\alpha y)=-f(J_\beta x,J_\beta y)=-f(J_\gamma x,J_\gamma y)=f(x,y),\quad \forall x, y \in
    V,
\end{equation}
where $(\alpha,\beta,\gamma)$ is a circular permutation of
$(1,2,3)$.
\end{defn}

Now, let us show the existence of the introduced bilinear forms on
$V$.

We denote $f\in \B_\alpha \subset \B(V)$ $(\alpha=1,2,3)$ when $f$
satisfies the conditions \eqref{16}.
Let us remark that $\B_\mathrm{H}(V)$ is a subspace of the vector
space $\B(V)$.
The projector $\Pi_H:\; \B(V) \rightarrow \B_\mathrm{H}(V)$ is
defined in \cite{AlMa} as follows
\begin{eqnarray}\label{19}
f \rightarrow (\Pi_\mathrm{H} f)(x,y):= \frac{1}{4} \left\{f(x,y)\right. &+& f(J_1x,J_1y) \nonumber \\
                                                    &+& \left. f(J_2x,J_2y)+f(J_3x,J_3y)\right\}.
\end{eqnarray}

For convenience we set $\Pi_0:=\Pi_\mathrm{H}$ and
$\B_0:=\B_\mathrm{H}(V)$. Clearly, $\Pi_0$ is a projector, \ie{}
$\Pi_0^2=\Pi_0$.

Analogously we define the operators: $\Pi_\alpha:\,
\B(V)\rightarrow \B_\alpha$, $\alpha=1,2,3$ as follows
\begin{eqnarray}\label{110}
f \rightarrow (\Pi_\alpha f)(x,y):=
\frac{1}{4}\left\{f(x,y)\right. &+& f(J_\alpha x,J_\alpha y) \nonumber \\
                                &-& \left.f(J_\beta x,J_\beta y)-f(J_\gamma x,J_\gamma y)\right\},
\end{eqnarray}
where $(\alpha,\beta,\gamma)$ is a circular permutation of
$(1,2,3)$.
It is not difficult to see that $\Pi_\alpha f \in \B_\alpha$,
$\alpha=1,2,3$.

In view of \eqref{19}--\eqref{110} the following proposition
holds:

\begin{prop}
The vector space $\B(V)$ admits the following decomposition
\[
\B(V)=\B_0 \oplus \B_1 \oplus \B_2 \oplus \B_3, \quad
\B_\alpha=\im\Pi_\alpha,\; \alpha=0,1,2,3,
\]
where the operators $\Pi_0, \Pi_1, \Pi_2$ and $ \Pi_3$ are
projectors with values in $\B(V)$ such that
\[
\begin{array} {ll}
\Pi_\alpha^2=\Pi_\alpha,\qquad & \Pi_0+\Pi_1+\Pi_2+\Pi_3=\Id,\\
\Pi_\alpha\circ\Pi_\beta=\Pi_\beta\circ\Pi_\alpha=0, \qquad &
\alpha\neq \beta;\quad \alpha,\beta\in \{0,1,2,3\}.
\end{array}
\]
\end{prop}

So, pseudo-Hermitian bilinear forms exist and moreover they are
three types in any vector space $V$ equipped with a hypercomplex
structure $H$, denoted by $(V,H)$.

Let $x$ determined by \eqref{11} and
\[
y=a^i\ddx+b^i\ddy+c^i\ddu+d^i\ddv, \quad i=1,2,\dots,n
\]
be arbitrary vectors on $V$. Following \cite{Wo}, we define as in
\cite{GrMaDi} a pseudo-Euclidean metric  of signature $(2n,2n)$ on
$V$ by a symmetric bilinear form $g$ as follows
\[
g(x,y):=\sum_{i=1}^n \left(-x^ia^i-y^ib^i+u^ic^i+v^id^i\right).
\]
Hence for the local basis $\left\{ \pd / \pd x^i, \pd / \pd y^i,
\pd / \pd u^i, \pd / \pd v^i \right\}$, $i=1,2$,$\dots,n$ on $V$
we have for $i,j\in\{1,2,\dots,n\}$
\begin{eqnarray}
-g\left(\frac{\pd}{\pd x^i},\frac{\pd}{\pd x^j}\right)&=&
-g\left(\frac{\pd}{\pd y^i},\frac{\pd}{\pd y^j}\right)=
g\left(\frac{\pd}{\pd u^i},\frac{\pd}{\pd u^j}\right)=
g\left(\frac{\pd}{\pd v^i},\frac{\pd}{\pd v^j}\right)=\delta_{ij},\nonumber \\[4pt]
g\left(\frac{\pd}{\pd x^i},\frac{\pd}{\pd y^j}\right)&=&
g\left(\frac{\pd}{\pd x^i},\frac{\pd}{\pd u^j}\right)=
g\left(\frac{\pd}{\pd x^i},\frac{\pd}{\pd v^j}\right)
=g\left(\frac{\pd}{\pd y^i},\frac{\pd}{\pd u^j}\right)\nonumber \\
&=&g\left(\frac{\pd}{\pd y^i},\frac{\pd}{\pd v^j}\right)
=g\left(\frac{\pd}{\pd u^i},\frac{\pd}{\pd v^j}\right)=0.\nonumber
\end{eqnarray}

Let us remark that if we denote $e_i=\pd / \pd x^i$
 $(i=1,2,\dots,n)$ then according to \eqref{12} the basis
\begin{equation}\label{116}
(e_1,e_2,\dots,e_n;J_1e_1,J_1e_2,\dots,J_1e_n;\dots;
    J_3e_1,J_3e_2,\dots,J_3e_n)
\end{equation}
is an {\em an admissible basis of $H$\/} and it is orthonormal
with respect to $g$.

Because of the properties
\begin{equation}\label{117}
g(J_1x,J_1y)=-g(J_2x,J_2y)=-g(J_3x,J_3y)=g(x,y),
\end{equation}
the pseudo-Euclidean metric $g$ is a symmetric pseudo-Hermit\-ian
bilinear form and $g \in \B_1$. Moreover,
$g_1(x,y):=g(J_1x,y)=-g(J_1y,x)$ coincides with the known K\"ahler
form with respect to $J_1$, \ie{} $\Phi(x,y):=g_1(x,y)$
\cite{GrMaDi}.

The associated bilinear forms $g_2(x,y):=g(J_2x,y)$ and
$g_3(x,y):=g(J_3x,y)$ of $g$ are symmetric and $\Phi \in \B_0, g
\in \B_1, g_2 \in \B_3, g_3 \in \B_2$, \ie{} the K\"ahler form
$\Phi$ is Hermitian and $g, g_2, g_3$ are pseudo-Hermitian of
different types, but they have the same signature $(2n,2n)$. Then
the structure $(H,G):=(H,g,\Phi,g_2,g_3)$ is called {\em a
hypercomplex pseudo-Hermitian structure\/} on $V$ \cite{GrMaDi}.

According to \cite{So}, the matrices that commute with $J_\alpha$
$(\alpha=1,2,3)$ are $A=\left(A_{ij}\right)$, $i,j\in
\{1,2,\dots,n\}$, where every $\left(A_{ij}\right)$ is a $(4\times
4)$-matrix of the form
\[
A_{ij}=\left(
 \begin{array}{cc}
            P &\; Q  \\
            -Q &\; P
 \end{array}\right),\quad
P=\left(
 \begin{array}{cc}
         a &\; b  \\
         -b &\; a
 \end{array}\right),\quad
Q=\left(
 \begin{array}{cc}
         c &\; d  \\
         d &\; -c
 \end{array}
 \right),\quad
a,b,c,d \in \mathbb{R}.
\]
The set of the $J_\alpha$-commuting matrices, that are invertible,
is a group which is isomorphic to $\mathrm{GL}(n,\mathbb{H})$.

The pseudo-Euclidean metric $g$ has a matrix with respect to the
basis \eqref{116} of the form $\mathbf{g}=\diag (g,g,\dots,g)$,
where
\[
g=\left( \begin{array}{cc}
     -I_2 & O_2  \\
     O_2 & I_2
\end{array} \right).
\]
The group preserving $\mathbf{g}$ is defined by the condition
$A^T\mathbf{g}A=\mathbf{g}$ for arbitrary $A\in
\mathrm{GL}(4n,\mathbb{R})$. It is clear that the group which
preserves $\mathbf{g}$ is $\mathrm{O}(2n,2n)$.

The structural group of $(V,H,G)$ has the property to preserve the
structures $J_\alpha$ and the metric $g$ (consequently $\Phi, g_2,
g_3$, too). Then this structural group is the intersection of
$\mathrm{GL}(n,\mathbb{H})$ and $\mathrm{O}(2n,2n)$.
 We get immediately that
\[
A\in \mathrm{GL}(n,\mathbb{H})\cap \mathrm{O}(2n,2n)
\quad\Leftrightarrow\quad a^2+b^2=1,\; c=d=0.
\]
Therefore $\mathrm{GL}(n,\mathbb{H})\cap \mathrm{O}(2n,2n)$ is an
1-parametrical group, \ie{} the elements $A_{ij}$ of $A$ depend on
1 real parameter.


\section{Almost $(H,G)$-structures on a manifold}

Let $(M,H)$ be an almost hypercomplex manifold \cite{AlMa}. We
suppose that $g$ is a symmetric tensor field of type $(0,2)$. If
it induces a pseudo-Hermitian inner product in $T_pM$, $p \in M$,
then $g$ is called {\em a pseudo-Hermitian metric on $M$\/}.
The structure $(H,G):=(J_1,J_2,J_3,g,\Phi,g_2,g_3)$ is called {\em
an almost hypercomplex pseudo-Hermitian structure on $M$\/} or in
short {\em an almost $(H,G)$-structure on $M$\/}.
The manifold $M$ equipped with $H$ and $G$, \ie{} $(M,H,G)$, is
called {\em an almost hypercomplex pseudo-Hermitian manifold\/},
or in short {\em an almost $(H,G)$-manifold\/}. The structural
tensors of the almost $(H,G)$-manifold are the three tensors of
type $(0,3)$ determined by
\begin{equation}\label{F}
F_\alpha (x,y,z)=g\bigl( \left( \nabla_x J_\alpha
\right)y,z\bigr)=\bigl(\nabla_x g_\alpha\bigr) \left( y,z
\right),\quad \alpha=1,2,3,
\end{equation}
where $\nabla$ is the Levi-Civita connection generated by $g$
\cite{GrMaDi}.

The properties of $H$ and $g$ imply the following properties of
$F_\alpha$:
\begin{equation}\label{32}
\begin{array}{l}
    F_1(x,y,z)=F_2(x,J_3y,z)+F_3(x,y,J_2z), \\
    F_2(x,y,z)=F_3(x,J_1y,z)+F_1(x,y,J_3z), \\
    F_3(x,y,z)=F_1(x,J_2y,z)-F_2(x,y,J_1z);
\end{array}
\end{equation}

\begin{equation}\label{34}
\begin{array}{l}
    F_1(x,y,z)=-F_1(x,z,y)=-F_1(x,J_1y,J_1z), \\
    F_2(x,y,z)=F_2(x,z,y)=F_2(x,J_2y,J_2z), \\
    F_3(x,y,z)=F_3(x,z,y)=F_3(x,J_3y,J_3z).
\end{array}
\end{equation}

Let us consider the Nijenhuis tensors $N_\alpha$ for $J_\alpha$
and $X,Y \in \X(M)$ given by \( N_\alpha(X,Y) = \left[J_\alpha
X,J_\alpha Y \right]
    -J_\alpha\left[J_\alpha X,Y \right]
    -J_\alpha\left[X,J_\alpha Y \right]
    -\left[X,Y \right].
\)
It is well known that the almost hypercomplex structure
$H=(J_\alpha)$ is a {\em hypercomplex structure\/} if $N_\alpha$
vanishes for each $\alpha=1,2,3$. Moreover, it is known that one
almost hypercomplex structure $H$ is hypercomplex if and only if
two of the structures $J_\alpha$ $(\alpha=1,2,3)$ are integrable.
This means that two of the tensors $N_\alpha$ vanish \cite{AlMa}.

Let us note that according to \eqref{117} the manifold $(M,J_1,g)$
is almost Hermitian and the manifolds $(M,J_\alpha,g)$,
$\alpha=2,3$, are almost complex manifolds with Norden metric (or
B-metric) \cite{GaBo,GaGrMi}. The basic classes of the mentioned
two types of manifolds for dimension $4n$ are:
\begin{equation}\label{cl-H}
\begin{array}{l}
\W_1(J_1):\; F_1(x,y,z)=-F_1(y,x,z), \\
\W_2(J_1):\; \mathop{\s}_{x,y,z}\bigl\{F_1(x,y,z)\bigr\}=0, \\
\W_3(J_1):\; F_1(x,y,z)=F_1(J_1x,J_1y,z),\quad \theta_1=0, \\
    \begin{array}{ll}
        \W_4(J_1):\; F_1(x,y,z)=&\frac{1}{2(2n-1)}
                \left\{g(x,y)\theta_1(z)-g(x,z)\theta_1(y)\right. \\

                &
                \left.-g(x,J_1y)\theta_1(J_1z)+g(x,J_1z)\theta_1(J_1y)
                \right\},
    \end{array}
\end{array}
\end{equation}
where $\theta_1(\cdot)=g^{ij}F_1(e_i,e_j,\cdot)$ for an arbitrary
basis $\{e_i\}_{i=1}^{4n}$ \cite{GrHe};
\begin{equation}\label{cl-N}
\begin{array}{l}
\begin{array}{ll}
\W_1(J_\alpha):\; F_\alpha(x,y,z)=&\frac{1}{4n}\bigl\{
g(x,y)\theta_\alpha(z)+g(x,z)\theta_\alpha(y)\bigr.\\
&\bigl.+g(x,J_\alpha y)\theta_\alpha(J_\alpha z)
    +g(x,J_\alpha z)\theta_\alpha(J_\alpha y)\bigr\},\\
\end{array}
\\
\W_2(J_\alpha):\; \mathop{\s}_{x,y,z}
\bigl\{F_\alpha(x,y,J_\alpha z)\bigr\}=0,\quad \theta_\alpha=0,\\
\W_3(J_\alpha):\; \mathop{\s}_{x,y,z}
\bigl\{F_\alpha(x,y,z)\bigr\}=0,
\end{array}
\end{equation}
where \(\theta_\alpha(z)=g^{ij}F_\alpha (e_i,e_j,z)\), \(\alpha
=2,3\), for an arbitrary basis \(\{e_i\}_{i=1}^{4n}\) and $\s $ is
the cyclic sum by three arguments \cite{GaBo}.

The special class $\W_0(J_\alpha):$ $F_\alpha=0$ $(\alpha =1,2,3)$
of the K\"ahler-type manifolds belongs to any other class within
the corresponding classification.

We say that an almost hypercomplex pseudo-Hermitian manifold is a
{\em pseudo-hyper-K\"ahler manifold\/}, if $\nabla J_\alpha=0\;
(\alpha=1,2,3)$ with respect to the Levi-Civita connection
generated by $g$ \cite{GrMaDi}.

Clearly, in this case we have $F_\alpha=0\, (\alpha=1,2,3)$ or the
manifold is K\"ahlerian with respect to $J_\alpha$, \ie{} $(M,H,G)
\in \W_0(J_\alpha)$.

Immediately we obtain
\begin{prop}\label{l33}
If $(M,H,G) \in \W_0(J_\alpha)\bigcap \W_0(J_\beta)$ then $(M,H,G)
\in \W_0(J_\gamma)$ for all cyclic permutations $(\alpha, \beta,
\gamma)$ of $(1,2,3)$ and $(M,H,G)$ is pseudo-hyper-K\"ahlerian.
$\hfill\Box$
\end{prop}

A basic property of the pseudo-hyper-K\"ahler manifolds is given
in \cite{GrMaDi} by the following
\begin{thm}[\cite{GrMaDi}]\label{R=0}
Each pseudo-hyper-K\"ahler manifold is a flat
pseudo-Rie\-mann\-ian manifold of signature $(2n,2n)$.
$\hfill\Box$
\end{thm}

As $g$ is an indefinite metric, there exist isotropic vector
fields $X$ on $M$, \ie{} \(g(X,X)=0\), \(X\neq 0\), \(X\in
\X(M)\). Following \cite{GRMa} we define the invariants
\begin{equation}\label{nJ}
\nJ{\alpha}= g^{ij}g^{kl}g\bigl( \left( \nabla_i J_\alpha \right)
e_k, \left( \nabla_j J_\alpha \right) e_l \bigr), \qquad
\alpha=1,2,3,
\end{equation}
where $\{e_i\}_{i=1}^{4n}$ is an arbitrary basis of $T_pM$, $p\in
M$. Let us remark that the invariant $\nJ{\alpha}$ is the scalar
square of the $(1,2)$-tensor $\nabla J_\alpha$.

\begin{defn}
We say that an $(H,G)$-manifold is:
\begin{enumerate}
\renewcommand{\labelenumi}{(\roman{enumi})}
\item isotropic K\"ahlerian with respect to $J_\alpha$
    if $\nJ{\alpha}=0$ for some $\alpha\in\{1,2,3\}$;
\item isotropic hyper-K\"ahlerian
    if it is isotropic K\"ahlerian with respect to every $J_\alpha$ of $H$.
\end{enumerate}
\end{defn}

Clearly, if $(M,H,G)$ is pseudo-hyper-K\"ahlerian, then it is an
isotropic hyper-K\"ahler manifold. The inverse statement does not
hold.


\section{A Lie group as a 4-dimensional $(H,G)$-mani\-fold}
\label{sec_3-2}

In \cite{GrMaMe1} is constructed an example of a 4-dimensional Lie
group equipped with a quasi-K\"ahler structure and Norden metric
$g$, \ie{} it is a $\W_3$-manifold according to \eqref{cl-N}.
There it is characterized with respect to $\nabla$ of $g$.

\begin{thm}[\cite{GrMaMe1}]\label{G}
Let $(\Lgr,J,g)$ be a 4-dimensional almost complex manifold with
Norden metric, where $\Lgr$ is a connected Lie group with a
corresponding Lie algebra $\mathfrak{l}$ determined by the global
basis of left invariant vector fields \\
$\{X_1,X_2,X_3,X_4\}$; $J$
is an almost complex structure defined by
\begin{equation}\label{J}
JX_1=X_3,\quad JX_2=X_4,\quad JX_3=-X_1,\quad JX_4=-X_2;
\end{equation}
$g$ is an invariant Norden metric determined by
\begin{equation}\label{g}
\begin{array}{l}
  g(X_1,X_1)=g(X_2,X_2)=-g(X_3,X_3)=-g(X_4,X_4)=1, \\[4pt]
  g(X_i,X_j)=0,\quad i\neq j; \qquad
  g\left([X_i,X_j],X_k\right)+g\left([X_i,X_k],X_j\right)=0.
\end{array}
\end{equation}
Then $(\Lgr,J,g)$ is a quasi-K\"ahler manifold with Norden metric
if and only if $\Lgr$ belongs to the 4-parametric family of Lie
groups determined by the conditions
\begin{equation}\label{[]4}
\begin{array}{ll}
[X_1,X_3]= \lambda_2 X_2+\lambda_4 X_4,\quad & [X_2,X_4]=\lambda_1
X_1+\lambda_3 X_3,
\\[4pt]
[X_2,X_3]= -\lambda_2 X_1-\lambda_3 X_4,\quad & [X_3,X_4]=-\lambda_4 X_1+\lambda_{3} X_2,\\[4pt]
[X_4,X_1]= \lambda_1 X_2+\lambda_4 X_3,\quad &
[X_2,X_1]=-\lambda_2 X_3+\lambda_1 X_4,
\end{array}
\end{equation}
where $\lambda_i\in \mathbb{R}$ $(i=1,2,3,4)$ and
$(\lambda_1,\lambda_2,\lambda_3,\lambda_4)\neq (0,0,0,0)$.
$\hfill\Box$
\end{thm}


The components of $\nabla$ are determined (\cite{GrMaMe1}) by
\eqref{[]4} and
\begin{equation}\label{invLC}
    \nabla_{X_i} X_j=\frac{1}{2}[X_i,X_j]\quad (i,j=1,2,3,4).
\end{equation}

Hence the components $R_{ijks}=R(X_i,X_j,X_k,X_s)$
$(i,j,k,s=1,2,3,4)$ of the curvature tensor $R$ on $(\Lgr,g)$ are:
\cite{GrMaMe1}
\begin{equation}\label{Rijks}
\begin{array}{ll}
    R_{1221}=-\frac{1}{4}\left(\lambda_1^2+\lambda_2^2\right),\qquad
    &
    R_{1331}=\frac{1}{4}\left(\lambda_2^2-\lambda_4^2\right),\\
    R_{1441}=-\frac{1}{4}\left(\lambda_1^2-\lambda_4^2\right),\qquad
    &
    R_{2332}=\frac{1}{4}\left(\lambda_2^2-\lambda_3^2\right),\\
    R_{2442}=\frac{1}{4}\left(\lambda_1^2-\lambda_3^2\right),\qquad
    &
    R_{3443}=\frac{1}{4}\left(\lambda_3^2+\lambda_4^2\right),\\
    R_{1341}=R_{2342}=-\frac{1}{4}\lambda_1\lambda_2,\qquad
    &
    R_{2132}=-R_{4134}=\frac{1}{4}\lambda_1\lambda_3,\\
    R_{1231}=-R_{4234}=\frac{1}{4}\lambda_1\lambda_4,\qquad
    &
    R_{2142}=-R_{3143}=\frac{1}{4}\lambda_2\lambda_3,\\
    R_{1241}=-R_{3243}=\frac{1}{4}\lambda_2\lambda_4,\qquad
    &
    R_{3123}=R_{4124}=\frac{1}{4}\lambda_3\lambda_4,\\
\end{array}
\end{equation}
and the scalar curvature $\tau$ on $(\Lgr,g)$ is \cite{GrMaMe1}
\begin{equation}\label{tau}
    \tau=-\frac{3}{2}\left(\lambda_1^2+\lambda_2^2-\lambda_3^2-\lambda_4^2\right).
\end{equation}

Now we introduce a hypercomplex structure $H=(J_1,J_2,J_3)$ by the
following way. At first, let $J_2$ be the given almost complex
structure $J$ by \eqref{J}. Secondly, we define an almost complex
structure $J_1$ as follows
\begin{equation}\label{J1}
\begin{array}{l}
J_1:\quad J_1X_1=X_2,\quad J_1X_2=-X_1,\quad J_1X_3=-X_4,\quad
J_1X_4=X_3.
\end{array}
\end{equation}
Finally, let the almost complex structure $J_3$ be the composition
of $J_1$ after $J_2$, \ie{} $J_3=J_1J_2$.

Then the introduced structure $(H,G)$ on $\Lgr$ has the properties
\eqref{13} and \eqref{117}. Hence we have the following
\begin{thm}\label{HG}
The manifold $(\Lgr,H,G)$ is an almost hypercomplex
pseudo-Her\-mitian
 manifold of dimension 4.$\hfill\Box$
\end{thm}

We continue by a characterization of the constructed manifold
$(\Lgr,H,G)$.

Let $(F_\alpha)_{ijk}=F_\alpha(X_i,X_j,X_k)$ and
$(\theta_\alpha)_i=\theta_\alpha (X_i)$ be the components of the
structural tensor $F_\alpha$ and its Lee form $\theta_\alpha$
($\alpha=1,2,3$), respectively. The nonzero components of $F_2$
are: \cite{GrMaMe1}

\begin{equation}\label{Fijk}
\begin{array}{l}
-(F_2)_{122}=-(F_2)_{144}=2(F_2)_{212}=2(F_2)_{221}=2(F_2)_{234}
\\
=2(F_2)_{243}=2(F_2)_{414}=-2(F_2)_{423}=-2(F_2)_{432}=2(F_2)_{441}=\lambda_1,
\\
2(F_2)_{112}=2(F_2)_{121}=2(F_2)_{134}=2(F_2)_{143}=-(F_2)_{211}
\\
=-(F_2)_{233}=-2(F_2)_{314}=2(F_2)_{323}=2(F_2)_{332}=-2(F_2)_{341}=\lambda_2,
\\
2(F_2)_{214}=-2(F_2)_{223}=-2(F_2)_{232}=2(F_2)_{241}=(F_2)_{322}
\\
=(F_2)_{344}=-2(F_2)_{412}=-2(F_2)_{421}=-2(F_2)_{434}=-2(F_2)_{443}=\lambda_3,
\\
-2(F_2)_{114}=2(F_2)_{123}=2(F_2)_{132}=-2(F_2)_{141}=-2(F_2)_{312}
\\
=-2(F_2)_{321}=-2(F_2)_{334}=-2(F_2)_{343}=(F_2)_{411}=(F_2)_{433}=\lambda_4.
\\
\end{array}
\end{equation}
Then we have $\theta_2=0$. By this way we confirm the statement in
\thmref{G} that the introduced manifold in \cite{GrMaMe1} is of
the basic class $\W_3$ with respect to $J_2$ within the
classification \eqref{cl-N}, \ie{}
\begin{equation}\label{W3-J2}
(\Lgr,J_2,g)\in \W_3(J_2).
\end{equation}

Having in mind  \eqref{g}--\eqref{invLC}, \eqref{J1} and
\eqref{F}, we obtain the nonzero components of $F_1$ as follows
\begin{equation}\label{F1}
    \begin{array}{l}
        (F_1)_{114}=-(F_1)_{123}=(F_1)_{132}=-(F_1)_{141}\\
        \phantom{(F_1)_{114}}
        =(F_1)_{213}=(F_1)_{224}=-(F_1)_{231}=-(F_1)_{242}=\frac{1}{2}\lambda_1;\\
        -(F_1)_{113}=-(F_1)_{124}=(F_1)_{131}=(F_1)_{142}\\
        \phantom{-(F_1)_{114}}
        =(F_1)_{214}=-(F_1)_{223}=(F_1)_{232}=-(F_1)_{241}=\frac{1}{2}\lambda_2;\\
        -(F_1)_{314}=(F_1)_{323}=-(F_1)_{332}=(F_1)_{341}\\
        \phantom{-(F_1)_{114}}
        =(F_1)_{413}=(F_1)_{424}=-(F_1)_{431}=-(F_1)_{442}=\frac{1}{2}\lambda_3;\\
        -(F_1)_{313}=-(F_1)_{324}=(F_1)_{331}=(F_1)_{342}\\
        \phantom{-(F_1)_{114}}
        =-(F_1)_{414}=(F_1)_{423}=-(F_1)_{432}=(F_1)_{441}=\frac{1}{2}\lambda_4.
    \end{array}
\end{equation}

Then we have
\[
    (\theta_1)_1=-\lambda_4,\quad
    (\theta_1)_2=\lambda_3,\quad
    (\theta_1)_3=-\lambda_2,\quad
    (\theta_1)_4=\lambda_1.
\]
Since $(\lambda_1,\lambda_2,\lambda_3,\lambda_4)\neq (0,0,0,0)$
then the 4-dimensional almost Hermitian manifold $(\Lgr,J_1,g)$ is
not K\"ahlerian and $\theta_1\neq 0$.

The validity of the property $F_1(X,Y,Z)=F_1(J_1X,J_1Y,Z)$ is
verified by us in virtue of \eqref{J1} and \eqref{F1}. It is
equivalent to the vanishing of the Nijenhuis tensor of $J_1$,
\ie{} $N_1=0$. According to \cite{GrHe} for dimension 4 we get
that the considered manifold belongs to the basic class
$\W_4(J_1)$ within \eqref{cl-H}, \ie{}
\begin{equation}\label{W4}
(\Lgr,J_1,g)\in \W_4(J_1).
\end{equation}

As it is known \cite{GrHe}, this class contains the conformally
K\"ahler manifolds of Hermitian type. The necessary and sufficient
condition a $\W_4(J_1)$-manifold to be locally or globally
conformally K\"ahlerian one is the Lee form $\theta_1$ to be
closed or exact. The basic components of $\dd{\theta_1}$ are:
\[
\begin{array}{ll}
\dd\theta_1(X_1,X_2)=\lambda_1^2+\lambda_2^2,\quad &
\dd\theta_1(X_2,X_4)=\dd\theta_1(X_3,X_1)=\lambda_1\lambda_4+\lambda_2\lambda_3,\\
\dd\theta_1(X_3,X_4)=-\lambda_3^2-\lambda_4^2,\quad  &
\dd\theta_1(X_1,X_4)=\dd\theta_1(X_2,X_3)=\lambda_1\lambda_3-\lambda_2\lambda_4.\\
\end{array}
\]

Hence, $\theta_1$ is not closed and therefore the constructed
$\W_4(J_1)$-manifold is not conformally K\"ahlerian.

Having in mind \eqref{32},   \eqref{Fijk},   \eqref{F1}, we
compute the following nonzero components of $F_3$:
\begin{equation}\label{F3}
\begin{array}{l}
(F_3)_{112}=(F_3)_{121}=-(F_3)_{134}=-(F_3)_{143}=-2(F_3)_{211}=-2(F_3)_{244}\\
\phantom{(F_3)_{112}}
=(F_3)_{413}=(F_3)_{431}=(F_3)_{424}=(F_3)_{442}=\frac{1}{2}\lambda_1,\\
2(F_3)_{122}=2(F_3)_{133}=-(F_3)_{212}=-(F_3)_{221}=(F_3)_{234}=(F_3)_{243}\\
\phantom{2(F_3)_{112}}
=-(F_3)_{313}=-(F_3)_{331}=-(F_3)_{324}=-(F_3)_{342}=\frac{1}{2}\lambda_2,\\
(F_3)_{213}=(F_3)_{231}=(F_3)_{224}=(F_3)_{242}=-(F_3)_{312}=-(F_3)_{321}\\
\phantom{(F_3)_{112}}
=(F_3)_{343}=(F_3)_{334}=-2(F_3)_{422}=-2(F_3)_{433}=\frac{1}{2}\lambda_3,\\
-(F_3)_{113}=-(F_3)_{124}=-(F_3)_{131}=-(F_3)_{142}=2(F_3)_{311}=2(F_3)_{344}\\
\phantom{-(F_3)_{112}}
=(F_3)_{412}=(F_3)_{421}=-(F_3)_{434}=-(F_3)_{443}=\frac{1}{2}\lambda_4.
\end{array}
\end{equation}

Hence, we establish directly that $\theta_3=0$ and
$\mathop{\s}_{i,j,k} (F_3)_{ijk}=0$. Therefore we obtain that the
considered manifold belongs to the basic class $\W_3(J_3)$, \ie{}
\begin{equation}\label{W3}
(\Lgr,J_3,g)\in \W_3(J_3).
\end{equation}

Let us summarize the conclusions \eqref{W3-J2}, \eqref{W4} and
\eqref{W3} in the following statement.
\begin{thm}
The constructed 4-dimensional almost hypercomplex
pseudo-Her\-mit\-ian manifold $(\Lgr,H,G)$ on the Lie group $\Lgr$
belongs to basic classes with respect to the three almost complex
structures of different types as follows
\[
(\Lgr,H,G) \in \W_4(J_1)\cap \W_3(J_2)\cap
\W_3(J_3).\qquad\qquad\qquad\qquad\qquad\quad\;\hfill\Box
\]
\end{thm}


The square norm $\nJ{\alpha}$ of $\nabla J_{\alpha}$ for an almost
complex structure $J_{\alpha}$ is defined in \cite{GRMa} by
\eqref{nJ}. Having in mind the definition
$F_{\alpha}(X,Y,Z)=g\bigl( \left( \nabla_X J_{\alpha}
\right)Y,Z\bigr)$ of the tensor $F_{\alpha}$,
we obtain the following equation for the square norm of $\nabla
J_{\alpha}$
\[
    \nJ{\alpha}=g^{ij}g^{kl}g^{pq}(F_{\alpha})_{ikp}(F_{\alpha})_{jlq},
\]
therefore
\begin{equation}\label{nJ=nF}
    \nJ{\alpha}=\norm{F_\alpha},\qquad \alpha=1,2,3.
\end{equation}

By virtue of \eqref{F1}, \eqref{Fijk}, \eqref{F3} we receive
immediately that
\[
    -2\nJ{1}=\nJ{2}=\nJ{3}=4\left(\lambda_1^2+\lambda_2^2-\lambda_3^2-\lambda_4^2\right).
\]

The last equations and Equation~\eqref{tau} imply
\begin{prop}\label{prop-iK}
\begin{enumerate}
\renewcommand{\labelenumi}{(\roman{enumi})}
\item If the manifold $(\Lgr,H,G)$ is isotropic K\"ahlerian with
respect to some $J_\alpha$ $(\alpha=1,2,3)$ then it is isotropic
hyper-K\"ahlerian; \item The manifold $(\Lgr,H,G)$ is isotropic
hyper-K\"ahlerian if and only if the condition
$\lambda_1^2+\lambda_2^2-\lambda_3^2-\lambda_4^2=0$
    holds;
\item The manifold $(\Lgr,H,G)$ is isotropic hyper-K\"ahlerian if
and only if it has zero scalar curvature $\tau$. $\hfill\Box$
\end{enumerate}
\end{prop}


The space of unitary invariants of order 2 for a 4-dimensional
Hermitian manifold is determined by the three quantities: $\tau$,
$\tau^*_1$ and $\norm{\nabla\Phi}=2\norm{\delta\Phi}$, where
$\tau^*_1=\frac{1}{2}g^{ij}g^{kl}R(X_i,J_1X_j,X_k,J_1X_l)$
\cite{GrHe}.

In other words, as $\norm{F_1}=\norm{\nabla\Phi}$ and
$\norm{\theta_1}=\norm{\delta\Phi}$ we get
\[
\begin{array}{l}
2\tau^*_1=-\norm{\theta_1}=\lambda_1^2+\lambda_2^2-\lambda_3^2-\lambda_4^2.
\end{array}
\]

Let us compute the associated scalar curvatures $\tau_\alpha^*$ on
$(\Lgr,J_\alpha,g)$ for $\alpha=2,3$ by
$\tau_\alpha^*:=g^{ij}g^{kl}R(X_i,X_k,J_\alpha X_l,X_j)$
\cite{GaGrMi}. Then, using \eqref{Rijks}, we obtain
\[
    \tau^*_2=\lambda_1\lambda_3+\lambda_2\lambda_4,\qquad \tau^*_3=\lambda_1\lambda_4-\lambda_2\lambda_3.
\]


Having in mind the definitions of the Nijenhuis tensors $N_\alpha$
of $J_\alpha$ ($\alpha=2,3$) and the commutators \eqref{[]4}, we
get the components $N_\alpha(X_i,X_j)$ and after that the square
norm of $N_\alpha$ as follows
\[
\norm{N_\alpha}=32\left(\lambda_1^2+\lambda_2^2-\lambda_3^2-\lambda_4^2\right),\quad
\alpha=2,3.
\]
It is clear, according to \propref{prop-iK} that the manifold
$(\Lgr,H,G)$ is isotropic hyper-K\"ahlerian and scalar flat if and
only if it has isotropic Nijenhuis tensors of $J_2$ and $J_3$.


\bigskip

\textit{Kostadin Gribachev, Mancho Manev\\
University of Plovdiv\\
Faculty of Mathematics and Informatics
\\
Department of Geometry\\
236 Bulgaria blvd.\\
Plovdiv 4003\\
Bulgaria}
\\
\texttt{e-mail: costas@uni-plovdiv.bg, mmanev@yahoo.com\\
http://www.fmi-plovdiv.org/manev}

\end{document}